\documentclass[letterpaper, 10 pt, conference]{ieeeconf}  

\IEEEoverridecommandlockouts                              

\usepackage{graphics} 
\usepackage{graphicx}
\usepackage{amsmath} 
\usepackage{amssymb}  
\usepackage{algorithm}
\usepackage{algorithmic}
\usepackage{color}
\usepackage{bm}
\usepackage{cite}

\title{\LARGE \bf
Information Theoretic Model Predictive Control on Jump Diffusion Processes
}

\author{Ziyi Wang$^{1}$, Grady Williams$^{1}$ and Evangelos A. Theodorou$^{1}$
\thanks{*This work was not supported by any organization}
\thanks{$^{1}$School of Aerospace Engineering, Georgia Tech, Atlanta, GA, USA}
}

\begin{document}

\maketitle
\thispagestyle{empty}
\pagestyle{empty}

\begin{abstract}

In this paper we present an information theoretic approach to stochastic optimal control problems for systems with compound Poisson noise. We generalize previous work on information theoretic path integral control to discontinuous dynamics with compound Poisson noise. We also derive a control update law of the same form using a stochastic optimization approach. We develop a sampling-based iterative model predictive control (MPC) algorithm. The proposed algorithm is parallelizable and when implemented on a Graphical Processing Unit (GPU) can run in real time. We test the performance of the proposed algorithm in simulation for two control tasks using a cartpole and a quadrotor system. Our simulations demonstrate improved performance of the new scheme and indicate the importance of incorporating the statistical characteristics of stochastic disturbances in the computation of the stochastic optimal control policies.

\end{abstract}

\section{INTRODUCTION}

Despite the maturity of the field of stochastic optimal control theory, the majority of the theoretical  and computational work  considers stochastic systems with Gaussian stochastic disturbances. This observation is  valid if one considers the lack of scalable and real time algorithms for control of high dimensional stochastic systems with disturbances that are far from being Gaussian and zero mean. Motivated by this lack of theory and algorithm on control of systems with more complex stochasticity, we consider dynamics with Gaussian and the more general compound Poisson noise.

Compound Poisson process, also known as the marked-jump process, is a doubly stochastic process where the stochasticity arises from both the jump time and amplitude \cite{hanson2007applied}. For simplicity, we use the term jump noise for compound Poisson noise. Processes with jumps have been widely used to describe the random evolution of, e.g., brain dynamics \cite{anvari2016disentangling}, of soil moisture dynamics \cite{daly2006probabilistic}, or of financial figures such as stock prices, market indices, and interest rates \cite{tankov2003financial}. For  application on dynamical systems,  jump stochastic  terms in the dynamics can  capture the discontinuities that arise due to phenomena such as gust or due to interactions of the system in consideration with the environment. Therefore, it is important that methods that deal with jump terms in the dynamics are developed.

The contributions of this paper are as follows:
\begin{itemize}
	\item We derive a novel algorithm for control of systems with jump noise from an information theoretic point of view. The new algorithm extends the capability of a previous scheme to handle a more complex form of stochasticity than the common Gaussian noise \cite{williams2016aggressive}.

	\item  We show the connection of the resulting scheme with an alternative approach to stochastic optimization that does not rely on importance sampling. With this equivalence established, convergence of the algorithm can be shown using techniques common in optimization literature \cite{zhou2014gradient}.

	\item We present an iterative MPC algorithm. The algorithm can utilize the parallel computing capabilities of the GPU, which means a large number of sampling trajectories can be propagated simultaneously and the algorithm can be implemented in real time. We implement the algorithm in simulation on a cartpole and quadrotor system with Gaussian and jump noise added, and we compare the performance of the algorithm against the path integral control based information theoretic MPC algorithm \cite{williams2016aggressive} that doesn't account for the jump noise.
\end{itemize}

The information theoretic approach we take in this paper is based on the path integral framework, which originated from Kappen and Theodorou's work \cite{kappen2005linear, theodorou2010generalized}. An iterative path integral control method, developed by Williams \cite{williams2016aggressive}, has been implemented for autonomous racing. This method uses the information theoretic notions of free energy and relative entropy, and obtains optimal control policy distribution through minimization of the Kullback-Leibler divergence (KL-Divergence) between a control induced probability measure and the optimal control policy induced probability measure. This approach allows for a solution to the stochastic optimal control problem using an importance sampling scheme.

The stochastic optimization approach is motivated by the lack of computational methods for stochastic optimal control. This approach is based on the stochastic approximation method \cite{borkar2009stochastic} where noisy observations are used to approximate stochastic functions. Optimization is then performed based on this noisy approximation to iteratively improve the solution. The derivation of this approach assumes a parameterized sampling distribution of entire trajectories based on system dynamics, reformulates the original problem with respect to the parameters, and obtains an update rule for these parameters through gradient descent.

The rest of this paper is organized as follows: in section II we provide the problem formulation. In section III we introduce the information theoretic approach to the stochastic optimal control problem. In section IV we introduce the stochastic optimization approach to the problem. Then we provide the MPC algorithm in section V. The simulation results are included in section VI. Finally, we conclude this paper in section VII.

\section{PROBLEM FORMULATION}

Consider a stochastic system with state $\mathbf{x}_t \in \mathbb{R}^n$ and control $\mathbf{u}_t \in \mathbb{R}^m$ at time t. We assume the dynamics also has additive noise from Brownian motion $\mathrm{d}\mathbf{w} \in \mathbb{R}^p$ and marked-jump process $\mathrm{d}\mathbf{P} \in\ \mathbb{R}^q$ with constant jump rate. We define $U \in \mathbb{R}^{m\times T}$ as the control sequence and $X \in \mathbb{R}^{n\times T}$ as the state trajectory over the time horizon $T$. We can formulate our stochastic optimal control problem as:

\begin{equation}
U^* = \arg\min_{U \in \mathcal{U}}\mathbb{E}_{\mathbb{Q}}\Big[\phi(\mathbf{x}_T,T)+\int^T_{t_0} \mathcal{L}(\mathbf{x}_t,\mathbf{u}_t,t)\mathrm{d}t \Big]
\end{equation}
where $\mathcal{U}$ is the set of admissible control sequences, and the expectation is taken with respect to the probability measure $\mathbb{Q}$ induced by the controlled dynamics:

\begin{equation}
\mathrm{d}\mathbf{x}_t = \mathbf{F}(\mathbf{x}_t,\mathbf{u}_t,t)\mathrm{d}t + \mathbf{B}(\mathbf{x}_t,t)\mathrm{d}\mathbf{w}^{(1)} + \mathbf{H}(\mathbf{x}_t,\mathcal{Q},t)\mathrm{d}\mathbf{P}^{(1)}
\end{equation}
with $\mathbb{E}[\mathrm{d}\mathbf{P}^{(1)}]=\nu^{(1)}\mathrm{d}t$ and $\nu^{(1)}$ is the jump rate. We assume zero mean normal distribution for the mark distribution, $\phi_{\mathcal{Q}}(q;t)\sim\mathcal{N}(0,\bm{\Sigma}_J)$. For the cost function we consider a state-dependent cost and a quadratic control cost:

\begin{equation}
\mathcal{L}(\mathbf{x}_t,\mathbf{u}_t,t) = q(\mathbf{x}_t,t) + \frac{1}{2} \mathbf{u}_t^\mathrm{T}\mathbf{R}(\mathbf{x}_t,t)\mathbf{u}_t
\end{equation}

We consider dynamics affine in control:

\begin{equation}
\mathbf{F}(\mathbf{x}_t,\mathbf{u}_t,t) = \mathbf{f}(\mathbf{x}_t,t) + \mathbf{G}(\mathbf{x}_t,t)\mathbf{u}_t
\end{equation}

The proof for existence and uniqueness of solution to the problem we are considering can be found in \cite{oksendal2005applied}.

\section{Information Theoretic Approach}

In this section we present the derivation of our sampling based stochastic trajectory optimization method for jump diffusion processes using an information theoretic approach.

Before we start the derivation we need to introduce two quantities from information theory that are the foundation to our derivation. First we define the \textit{Free Energy} of a system as:

\begin{equation}
\mathcal{F}(S(X)) = -\log\Big(\mathbb{E}_\mathbb{P}\Big[\exp\big(-\frac{1}{\lambda}S(X)\big)\Big]\Big)
\end{equation}
where $\lambda \in \mathbb{R}^+$ is called the inverse temperature, and $S(X)$ is the state-dependent cost of a trajectory, $S(X) = \phi(\mathbf{x}_T,T) + \int^T_{t_0} q(\mathbf{x}_t,t)\mathrm{d}t$. The expectation is taken with respect to $\mathbb{P}$, which is the probability measure induced by the uncontrolled dynamics:

\begin{equation}
\mathrm{d}\mathbf{x}_t = \mathbf{f}(\mathbf{x}_t,t)\mathrm{d}t + \mathbf{B}(\mathbf{x}_t,t)\mathrm{d}\mathbf{w}^{(0)} + \mathbf{H}(\mathbf{x}_t,\mathcal{Q},t)\mathrm{d}\mathbf{P}^{(0)}
\end{equation}
with $\mathbb{E}[\mathrm{d}\mathbb{P}^{(0)}] = \nu^{(0)}\mathrm{d}t$ and the same mark distribution as the controlled dynamics. Next let $\mathbb{M},\mathbb{N}$ be two probability distributions that are absolutely continuous with each other. Then the \textit{KL-Divergence} between them is:

\begin{equation}
\mathbb{D}_{KL}(\mathbb{M}\parallel\mathbb{N}) = \mathbb{E}_\mathbb{M}\Big[\log\Big(\frac{\mathrm{d}\mathbb{M}}{\mathrm{d}\mathbb{N}}\Big)\Big]
\end{equation}

The KL-Divergence provides a measure of how one probability distribution diverges from a second and can be roughly thought of as the distance between two probability distributions, although it is not symmetric. The KL-Divergence is useful for defining optimization objectives.

Now suppose probability distributions $\mathbb{Q}$ and $\mathbb{P}$ as defined previously are absolutely continuous with each other, we can make the following observation:

\begin{equation}
\begin{split}
\mathcal{F}(S(X)) &= -\log\Big(\mathbb{E}_\mathbb{P}\Big[\exp\big(-\frac{1}{\lambda}S(X)\big)\Big]\Big) \\
&= -\log\Big(\mathbb{E}_{\mathbb{Q}}\Big[\exp\big(-\frac{1}{\lambda}S(X)\big)\frac{\mathrm{d}\mathbb{P}}{\mathrm{d}\mathbb{Q}} \Big] \Big)
\end{split}
\end{equation}
where we changed the expectation by multiplying by $1=\frac{\mathrm{d}\mathbb{Q}}{\mathrm{d}\mathbb{Q}}$. Since the negative logarithm is a convex function, we can apply Jensen's inequality and obtain:

\begin{equation}
\mathcal{F}(S(X)) \leq -\mathbb{E}_{\mathbb{Q}}\Big[\log\Big(\exp\big(-\frac{1}{\lambda}S(X)\big)\frac{\mathrm{d}\mathbb{P}}{\mathrm{d}\mathbb{Q}} \Big) \Big]
\end{equation}

The right hand side can be simplified as:

\begin{equation}
\begin{split}
RHS&= \frac{1}{\lambda}\mathbb{E}_{\mathbb{Q}}\Big[S(X) +  \lambda\log\Big(\frac{\mathrm{d}\mathbb{P}}{\mathrm{d}\mathbb{Q}} \Big) \Big] \\
&= \frac{1}{\lambda} \Big(\mathbb{E}_{\mathbb{Q}}[S(X)] + \lambda \mathbb{D}_{KL}(\mathbb{Q}\parallel\mathbb{P}) \Big)
\end{split}
\end{equation}

Substituting the terms back to (9):

\begin{equation}
 \lambda\mathcal{F}(S(X)) \leq \mathbb{E}_{\mathbb{Q}}[S(X)] + \lambda \mathbb{D}_{KL}(\mathbb{Q}\parallel\mathbb{P})
\end{equation}

To find the KL-Divergence between $\mathbb{Q}$ and $\mathbb{P}$, we need $\frac{\mathrm{d}\mathbb{Q}}{\mathrm{d}\mathbb{P}}$, which can be found using Girsanov's theorem \cite{hanson2007applied}:

\begin{equation}
\begin{split}
\frac{\mathrm{d}\mathbb{Q}}{\mathrm{d}\mathbb{P}} &= \exp\Big(\frac{1}{2}\int^T_{t_0}\mathbf{u}_t^\mathrm{T}\mathbf{G}(\mathbf{x}_t,t)^\mathrm{T}\Sigma(\mathbf{x}_t,t)^{-1}\mathbf{G}(\mathbf{x}_t,t)\mathbf{u}_t\mathrm{d}t \\
&+ \int^T_{t_0}\mathbf{u}_t^\mathrm{T}\mathbf{G}(\mathbf{x}_t,t)^\mathrm{T}\Sigma(\mathbf{x}_t,t)^{-1}\mathbf{B}(\mathbf{x}_t,t)\mathrm{d}\mathbf{w}^{(1)}\\
&-\int^T_{t_0}((\gamma^J-1)\nu^{(0)})\mathrm{d}t\Big)\cdot\prod_{k=1}^{\mathbf{P}^{(0)}(t)}\gamma^J(T_k^-)\gamma^M(\mathcal{Q}_k,T_k^-)
\end{split}
\end{equation}
where $\Sigma(\mathbf{x}_t,t) = \mathbf{B}(\mathbf{x}_t,t)\mathbf{B}(\mathbf{x}_t,t)^T$, $\gamma^J(t)$ is the ratio of jump rates in the two dynamics, $\int_0^T\nu^{(1)}\mathrm{d}t=\int_0^T\gamma^J(t)\nu^{(0)}\mathrm{d}t$, and $\gamma^M(q;t)$ is the scaling between the mark distributions, $\int_{\mathcal{Q}_1}\phi_{\mathcal{Q}}^{(1)}(q;t)\mathrm{d}q=\int_{\mathcal{Q}_0}\gamma^M(q;t)\phi_{\mathcal{Q}}^{(0)}(q;t)\mathrm{d}q=1$.

Here we consider the case where the change of measure only includes changes in drift, and the jump rates and mark distributions are the same. Therefore, both $\gamma^J$ and $\gamma^M$ have the value 1, and the last two terms can be dropped. Additionally, since $\mathrm{d}\mathbf{w}^{(1)}$ is a Brownian motion with respect to $\mathbb{Q}$, we get $\mathbb{E}_{\mathbb{Q}}\Big[\int_0^T\mathrm{d}\mathbf{w}^{(1)}\Big]=0$. The KL-Divergence then simplifies to:

\begin{equation}
\begin{split}
&\mathbb{D}_{KL}(\mathbb{Q}\parallel\mathbb{P}) =\\
&\mathbb{E}_\mathbb{Q}\Big[\frac{1}{2}\int^T_{t_0}\mathbf{u}_t^\mathrm{T}\mathbf{G}(\mathbf{x}_t,t)^\mathrm{T}\Sigma(\mathbf{x}_t,t)^{-1}\mathbf{G}(\mathbf{x}_t,t)\mathbf{u}_t\mathrm{d}t\Big]
\end{split}
\end{equation}

Using this result, if we assume the control cost matrix has the form:

\begin{equation}
\mathbf{R}(\mathbf{x}_t,t) = \lambda\mathbf{G}(\mathbf{x}_t,t)^\mathrm{T}\Sigma(\mathbf{x}_t,t)^{-1}\mathbf{G}(\mathbf{x}_t,t)
\end{equation}
we get the following form on the right hand side of equation (11):

\begin{equation}
\begin{split}
\mathbb{E}_\mathbb{Q}[S(X)]&+\lambda \mathbb{D}_{KL}(\mathbb{Q}\parallel\mathbb{P}) =\\
&\mathbb{E}_\mathbb{Q}\Big[S(X)+\frac{1}{2}\int^T_{t_0}\mathbf{u}_t^\mathrm{T}\mathbf{R}(\mathbf{x}_t,t)\mathbf{u}_t\mathrm{d}t\Big]
\end{split}
\end{equation}

Note that this is equivalent to the cost function in (1). With this we have shown that the free energy serves as the lower bound for our stochastic optimal control problem, and we can rewrite (11) as a minimization problem:

\begin{equation}
 \lambda\mathcal{F}(S(X)) = \inf_\mathbb{Q}\Big[\mathbb{E}_\mathbb{Q}[S(X)]+\lambda \mathbb{D}_{KL}(\mathbb{Q}\parallel\mathbb{P})\Big]
\end{equation}

In this minimization problem we have a state cost and a control cost in the form of KL-Divergence, which penalizes deviation from the uncontrolled distribution. We now define the optimal measure that achieves the lower bound as:

\begin{equation}
\frac{\mathrm{d}\mathbb{Q}^*}{\mathrm{d}\mathbb{P}} = \frac{\exp(-\frac{1}{\lambda}S(X))}{\mathbb{E}_\mathbb{P}[\exp(-\frac{1}{\lambda}S(X))]}
\end{equation}

This result can be easily verified by plugging it into (11) and is derived in [9]. With this we can solve the minimization problem defined by (16) by moving the probability distribution $\mathbb{Q}$ induced by some control as close to the optimal distribution as possible. The distance can be represented by the KL-Divergence between the two distributions and the problem becomes:

\begin{equation}
U^* = \arg\min_{U\in\mathcal{U}} \mathbb{D}_{KL}(\mathbb{Q}^*\parallel\mathbb{Q})
\end{equation}

\subsection{KL-Divergence Minimization}
Applying the definition of KL-Divergence we have:

\begin{equation}
\begin{split}
\mathbb{D}_{KL}(\mathbb{Q}^*\parallel\mathbb{Q}) &= \mathbb{E}_{\mathbb{Q}^*}\Big[\log\Big(\frac{\mathrm{d}\mathbb{Q}^*}{\mathrm{d}\mathbb{Q}}\Big)\Big] \\
&= \mathbb{E}_{\mathbb{Q}^*}\Big[\log\Big(\frac{\mathrm{d}\mathbb{Q}^*}{\mathrm{d}\mathbb{P}}\frac{\mathrm{d}\mathbb{P}}{\mathrm{d}\mathbb{Q}}\Big)\Big]
\end{split}
\end{equation}

We already have $\frac{\mathrm{d}\mathbb{Q}^*}{\mathrm{d}\mathbb{P}}$ from its definition. For $\frac{\mathrm{d}\mathbb{P}}{\mathrm{d}\mathbb{Q}}$, we can use Girsanov's theorem:

\begin{equation}
\begin{split}
\frac{\mathrm{d}\mathbb{P}}{\mathrm{d}\mathbb{Q}} &= \exp\Big(\frac{1}{2}\int^T_{t_0}\mathbf{u}_t^\mathrm{T}\mathbf{G}(\mathbf{x}_t,t)^\mathrm{T}\Sigma(\mathbf{x}_t,t)^{-1}\mathbf{G}(\mathbf{x}_t,t)\mathbf{u}_t\mathrm{d}t \\
&- \int^T_{t_0}\mathbf{u}_t^\mathrm{T}\mathbf{G}(\mathbf{x}_t,t)^\mathrm{T}\Sigma(\mathbf{x}_t,t)^{-1}\mathbf{B}(\mathbf{x}_t,t)\mathrm{d}\mathbf{w}^{(0)}\Big)
\end{split}
\end{equation}

Setting the terms inside the exponential as $\mathcal{D}(X,U)$ and plugging the results back in (19) we have:

\begin{equation}
\begin{split}
&\mathbb{D}_{KL}(\mathbb{Q}^*\parallel\mathbb{Q}) =\\
& \mathbb{E}_{\mathbb{Q}^*}\Big[-\frac{1}{\lambda}S(X)-\log(\mathbb{E}_\mathbb{P}[\exp(-\frac{1}{\lambda}S(X))])+\mathcal{D}(X,U)\Big]
\end{split}
\end{equation}

Since $S(X)$ is not dependent on the control we can drop the first two terms from the minimization. Now we discretize the control as step functions $\mathbf{u}_t=\mathbf{u}_j$ if $j\Delta t\leq t < (j+1)\Delta t$ with $j=\{0,1,\cdots,N-1\}$. Then we have:

\begin{equation}
\begin{split}
&\mathcal{D}(X,U) =\\
&\sum_{j=0}^{N-1}\Bigg(\frac{1}{2}\mathbf{u}_j^\mathrm{T}\int^{t_{j+1}}_{t_j} \mathcal{G}(\mathbf{x}_t,t)\mathrm{d}t \mathbf{u}_j - \mathbf{u}_j^\mathrm{T}\int^{t_{j+1}}_{t_j} \mathcal{B}(\mathbf{x}_t,t) \mathrm{d}\mathbf{w}^{(0)}\Bigg)
\end{split}
\end{equation}
where

\begin{equation}
\mathcal{G}(\mathbf{x}_t,t) = \mathbf{G}(\mathbf{x}_t,t)^\mathrm{T}\Sigma(\mathbf{x}_t,t)^{-1}\mathbf{G}(\mathbf{x}_t,t)
\end{equation}
\begin{equation}
\mathcal{B}(\mathbf{x}_t,t) = \mathbf{G}(\mathbf{x}_t,t)^\mathrm{T}\Sigma(\mathbf{x}_t,t)^{-1}\mathbf{B}(\mathbf{x}_t,t)
\end{equation}
\begin{equation}
N = T/\Delta t
\end{equation}

Note that each $\mathbf{u}_j$ does not depend on the trajectory taken, so we can taken them out of the expectation:

\begin{equation}
\begin{split}
\mathbb{E}_{\mathbb{Q}^*}\Big[\mathcal{D}(X,U)\Big] &= \sum_{j=0}^{N-1}\Bigg(\frac{1}{2}\mathbf{u}_j^\mathrm{T}\mathbb{E}_{\mathbb{Q}^*}\Big[\int^{t_{j+1}}_{t_j} \mathcal{G}(\mathbf{x}_t,t)\mathrm{d}t\Big] \mathbf{u}_j \\
&- \mathbf{u}_j^\mathrm{T}\mathbb{E}_{\mathbb{Q}^*}[\int^{t_{j+1}}_{t_j} \mathcal{B}(\mathbf{x}_t,t) \mathrm{d}\mathbf{w}^{(0)}]\Bigg)
\end{split}
\end{equation}

We can approximate the two integrals for small enough $\Delta t$ as:

\begin{equation}
\int^{t_{j+1}}_{t_j} \mathcal{G}(\mathbf{x}_t,t)\mathrm{d}t \approx \mathcal{G}(\mathbf{x}_{t_j},t_j) \Delta t
\end{equation}
\begin{equation}
\int^{t_{j+1}}_{t_j} \mathcal{B}(\mathbf{x}_t,t) \mathrm{d}\mathbf{w}^{(0)} \approx \mathcal{B}(\mathbf{x}_{t_j},t_j) \epsilon^{(0)}_j\sqrt{\Delta t}
\end{equation}
where $\epsilon^{(0)}_j$ is a vector with standard normal variable in each entry, $\epsilon^{(0)}_j\sim\mathcal{N}(0,\bm{\Sigma}_D)$. Then we can find $\mathbf{u}_j^*$ by taking the gradient with respect to $\mathbf{u}_j$, setting it to zero and solving for $\mathbf{u}_j$. The optimal control is found as:

\begin{equation}
\mathbf{u}^*_j = \frac{1}{\Delta t} \mathbb{E}_{\mathbb{Q}^*}\Big[\mathcal{G}(\mathbf{x}_{t_j},t_j)\Big]^{-1}\mathbb{E}_{\mathbb{Q}^*}\Big[\mathcal{B}(\mathbf{x}_{t_j},t_j)\epsilon^{(0)}_j\sqrt{\Delta t}\Big]
\end{equation}

\subsection{Importance Sampling}
We have obtained the optimal control in the form of expectation with respect to the optimal distribution. We can't sample from the optimal distribution, but we can sample from the uncontrolled distribution $\mathbb{P}$ to approximate the controls. Therefore, we need to change the expectation through multiplying by $\frac{\mathrm{d}\mathbb{P}}{\mathrm{d}\mathbb{P}}$ and using the Radon-Nikodym derivative $\frac{\mathrm{d}\mathbb{Q}^*}{\mathrm{d}\mathbb{P}}$:

\begin{equation}
\begin{split}
\mathbf{u}^*_j &= \frac{1}{\Delta t} \mathbb{E}_\mathbb{P}\Big[\frac{\exp(-\frac{1}{\lambda}S(X))\mathcal{G}(\mathbf{x}_{t_j},t_j)}{\mathbb{E}_\mathbb{P}[\exp(-\frac{1}{\lambda}S(X))]}\Big]^{-1}\\
&\cdot\mathbb{E}_\mathbb{P}\Big[\frac{\exp(-\frac{1}{\lambda}S(X)) \mathcal{B}(\mathbf{x}_{t_j},t_j)\epsilon^{(0)}_j\sqrt{\Delta t}}{\mathbb{E}_\mathbb{P}[\exp(-\frac{1}{\lambda}S(X))]}\Big]
\end{split}
\end{equation}

The equation can be further simplified since $\mathcal{G}(\mathbf{x}_{t_j},t_j)$ and $\mathcal{B}(\mathbf{x}_{t_j},t_j)$ are deterministic at time $t_j$:

\begin{equation}
\begin{split}
&\mathbf{u}^*_j =\\
&\frac{1}{\Delta t}\mathcal{G}(\mathbf{x}_{t_j},t_j)^{-1} \mathcal{B}(\mathbf{x}_{t_j},t_j) \mathbb{E}_\mathbb{P}\Big[\frac{\exp(-\frac{1}{\lambda}S(X)) \epsilon^{(0)}_j\sqrt{\Delta t}}{\mathbb{E}_\mathbb{P}[\exp(-\frac{1}{\lambda}S(X))]}\Big]
\end{split}
\end{equation}

Note that the expectations are taken with respect to the uncontrolled dynamics. This is not ideal since it means waiting for random Gaussian and jump noise to generate a meaningful trajectory. Therefore, we need to change the sampling distribution to the control induced distribution. In addition, we can also change the sampling variance to $\bm{\Sigma}_D^{(1)} = c \bm{\Sigma}_D^{(0)}$ to increase the state space explored. To perform importance sampling we multiply by $\frac{\mathrm{d}\mathbb{Q}}{\mathrm{d}\mathbb{Q}}$ and change from the zero mean $\epsilon^{(0)}_j \sqrt{\Delta t}$ to the non zero mean $\mathbf{G}(\mathbf{x})\mathbf{u}_j\Delta t + \epsilon^{(1)}_j \sqrt{\Delta t}$:

\begin{equation}
\begin{split}
\mathbf{u}^*_j &= \mathbf{u}_j + \frac{1}{\Delta t}\mathcal{G}(\mathbf{x}_{t_j},t_j)^{-1} \mathcal{B}(\mathbf{x}_{t_j},t_j)\\
&\cdot \mathbb{E}_\mathbb{Q}\Big[\frac{\exp(-\frac{1}{\lambda}S(X)) \epsilon^{(1)}_j\sqrt{\Delta t}\frac{\mathrm{d}\mathbb{P}}{\mathrm{d}\mathbb{Q}}}{\mathbb{E}_\mathbb{Q}[\exp(-\frac{1}{\lambda}S(X))\frac{\mathrm{d}\mathbb{P}}{\mathrm{d}\mathbb{Q}}]}\Big]
\end{split}
\end{equation}

We can use Girsanov's theorem again to get $\frac{\mathrm{d}\mathbb{P}}{\mathrm{d}\mathbb{Q}}$:

\begin{equation}
\begin{split}
\frac{\mathrm{d}\mathbb{P}}{\mathrm{d}\mathbb{Q}} &= \exp\Bigg(-\frac{1}{2} \sum_{j=0}^{N-1}\Big(\mathbf{u}_j^\mathrm{T}\mathbf{G}(\mathbf{x}_{t_j})^\mathrm{T}\Sigma^{-1}\mathbf{G}(\mathbf{x}_{t_j})\mathbf{u}_j \Delta t \\
&+ \mathbf{u}_j^\mathrm{T}\mathbf{G}(\mathbf{x}_{t_j})^\mathrm{T}\Sigma^{-1}\mathbf{B}(\mathbf{x}_{t_j})\epsilon_j^{(1)}\sqrt{\Delta t}\\
&+ (1-c^{-1})\epsilon_j^{(1)\mathrm{T}}\mathbf{B}(\mathbf{x}_{t_j})^\mathrm{T}\Sigma^{-1}\mathbf{B}(\mathbf{x}_{t_j})\epsilon_j^{(1)} \Delta t \Big)\Bigg)
\end{split}
\end{equation}

The last terms comes from the change of sampling variance and the detailed derivation can be found in \cite{williams2015model}. The addition of these terms can be added into the state cost:

\begin{equation}
\tilde{S}(X) = \phi(\mathbf{x}_{t_N}, t_N) + \sum_{j=0}^{N-1} \tilde{q}(\mathbf{x}_{t_j},\mathbf{u}_j,t_j) \Delta t
\end{equation}

\begin{equation}
\begin{split}
&\tilde{q}(\mathbf{x}_{t_j},\mathbf{u}_j,t_j) = q(\mathbf{x}_{t_j},t_j) + \frac{1}{2} \mathbf{u}_j^\mathrm{T}\mathbf{R}\mathbf{u}_j + \lambda \mathbf{u}_j^\mathrm{T}\mathcal{B}\frac{\epsilon_j^{(1)}}{\sqrt{\Delta t}} \\
&+ \frac{1}{2} \lambda (1-c^{-1})\epsilon_j^{(1)\mathrm{T}}\mathbf{B}(\mathbf{x}_{t_j})^\mathrm{T}\Sigma^{-1}\mathbf{B}(\mathbf{x}_{t_j})\epsilon_j^{(1)}/\Delta t
\end{split}
\end{equation}

With the new state cost we can obtain the final expression of optimal control update rule:

\begin{equation}
\mathbf{u}^*_j = \mathbf{u}_j + \mathcal{G}(\mathbf{x}_{t_j},t_j)^{-1}\mathcal{B}(\mathbf{x}_{t_j},t_j) \Big(\frac{\mathbb{E}_\mathbb{Q}[\exp(-\frac{1}{\lambda}\tilde{S}(X)) \frac{\epsilon^{(1)}_j}{\sqrt{\Delta t}}]}{\mathbb{E}_\mathbb{Q}[\exp(-\frac{1}{\lambda}\tilde{S}(X))]} \Big)
\end{equation}

The term inside the square brackets is approximated as:

\begin{equation}
\frac{\sum_{m=1}^M\exp(-\frac{1}{\lambda}\tilde{S}(X^m)) \frac{\epsilon_j^m}{\sqrt{\Delta t}}}{\sum_{m=1}^M\exp(-\frac{1}{\lambda}\tilde{S}(X^m))}
\end{equation}
using $M$ sample trajectories.

\section{Stochastic Optimization Approach}
\subsection{Problem Reformulation}
Consider a system with the same definition of state, control and dynamics as in the previous section, the optimal control problem can be defined in the same way:

\begin{equation}
U^* = \arg\min_{U\in\mathcal{U}}\mathbb{E}[J(X,U)]\label{eq:1}
\end{equation}
where $J:\mathbb{R}^{n \times T} \times \mathbb{R}^{m \times T}\rightarrow\mathbb{R}$ is an arbitrary cost function. We can introduce an exponential shape function $L(y)=\exp(y)$ to redefine the optimal control problem as a maximization problem:

\begin{equation}
U^* = \arg\max_{U\in\mathcal{U}}\mathbb{E}\Big[L\Big(-\frac{1}{\lambda} J(X,U)\Big)\Big]
\end{equation}

The expectation is taken over the control policy, which is parameterized by a set of parameters $\theta\in\Theta$ that we have control over. Finally, since $\ln:\mathbb{R}^+\rightarrow\mathbb{R}$ is a strictly increasing function, reformulating the maximization problem does not change the solution to the original optimal control problem:

\begin{equation}
\begin{split}
\theta^* &= \arg\max_{\theta\in\Theta}\ln\Big(\mathbb{E}\Big[L\Big(-\frac{1}{\lambda} J(X,U)\Big)\Big]\Big)\\
&= \arg\max_{\theta\in\Theta} l(\theta)
\end{split}
\end{equation}

\subsection{Probability Distribution Parameterization}

Assume a time discretization of control policy with step functions, and the stochastic control policy at each time instant has additive Gaussian and jump noise around some mean, $\mathbf{u}_j=\bm{\mu}_j+\epsilon_{D,j}+\epsilon_{J,j}\Delta\mathbf{P}$. We have the diffusion term $\epsilon_{D,j}\sim\mathcal{N}(0,\bm{\Sigma}_D)$ and the jump term $\epsilon_{J,j}\sim\mathcal{N}(0,\bm{\Sigma}_J)$ with known variances, and $\mathbb{E}[\Delta P]=\nu\Delta t$ with known jump rate $\nu$. Assuming stochasticity enters the system through the control channels, the probability density/mass function of each trajectory can be expressed as $p(X,U;\theta)=p(U;\theta)$ since the dynamics is deterministic. Assume the noise at each time instant is i.i.d., then the pdf/pmf of each trajectory can be expressed as:

\begin{equation}
p(U;\theta) = \prod^{N-1}_{j=0} p(\mathbf{u}_j;\theta_j)
\end{equation}

Although there is no closed form expression for the pdf/pmf of the entire trajectory, the pdf/pmf at each time instant can be written out explicitly. For small enough $\Delta t$ such that $\nu\Delta t \ll 1$, the zero-one jump law \cite{hanson2007applied} applies and the jump noise has Bernoulli distribution with probability of jump being $\nu\Delta t$. In addition, when jump occurs, $\mathbf{u}_j\sim\mathcal{N}(\bm{\mu}_j,\bm{\Sigma}_D+\bm{\Sigma}_J)$ is normally distributed with the variance as the sum of diffusion and jump noise since the sum of two normally distributed random variable is still normal. Therefore the pdf/pmf can be expressed as:

\begin{equation}
\begin{split}
p(\mathbf{u}_j;\theta_j) &= I_j(\nu\Delta t)\Bigg(\frac{1}{\sqrt{(2\pi)^n|\bm{\Sigma_D}+\bm{\Sigma_J}|}}\\
&\exp\Big(-\frac{1}{2}(\mathbf{u}_j-\bm{\mu}_j)^\mathrm{T}(\bm{\Sigma}_D+\bm{\Sigma}_J)^{-1}(\mathbf{u}_j-\bm{\mu}_j)\Big)\Bigg)\\
&+ (1-I_j)(1-\nu\Delta t)\Bigg(\frac{1}{\sqrt{(2\pi)^n|\bm{\Sigma_D}|}}\\
&\exp\Big(-\frac{1}{2}(\mathbf{u}_j-\bm{\mu}_j)^\mathrm{T}\bm{\Sigma}_D^{-1}(\mathbf{u}_j-\bm{\mu}_j)\Big)\Bigg)\\
&= h(\mathbf{u}_j)\exp\Big(\theta_j^\mathrm{T}T(\mathbf{u}_j)-A(\theta_j)\Big)
\end{split}
\end{equation}
where
\begin{equation}
\begin{split}
h(\mathbf{u}_j) &= I_j(\nu\Delta t)\Bigg(\frac{1}{\sqrt{(2\pi)^n|\bm{\Sigma}_D+\bm{\Sigma}_J|}}\\
&\exp\Big(-\frac{1}{2}\mathbf{u}_j^\mathrm{T}(\bm{\Sigma}_D+\bm{\Sigma}_J)^{-1}\mathbf{u}_j\Big)\Bigg) \\
&+ (1-I_j)(1-\nu\Delta t)\Bigg(\frac{1}{\sqrt{(2\pi)^n|\bm{\Sigma}_D|}}\\
&\exp\Big(-\frac{1}{2}\mathbf{u}_j^\mathrm{T}\bm{\Sigma}_D^{-1}\mathbf{u}_j\Big)\Bigg)\\
A(\theta_j) &= \frac{1}{2} \bm{\mu}_j^\mathrm{T}(\bm{\Sigma}_D+I_j\bm{\Sigma}_J)^{-1}\bm{\mu}_j \\
T(\mathbf{u}_j) &= (\bm{\Sigma}_D+I_j\bm{\Sigma}_J)^{-\frac{1}{2}} \mathbf{u}_j \\
\theta_j &= (\bm{\Sigma}_D+I_j\bm{\Sigma}_J)^{-\frac{1}{2}}\bm{\mu}_j
\end{split}
\end{equation}

The term $I_j$ term is an indicator function with $I_j=1$ when jump occurs and $I_j=0$ when there is no jump.

\subsection{Gradient Descent}
Now the gradient can be taken for the step direction at each time step j:
\begin{equation}
\begin{split}
\nabla_{\theta_j} l(\theta) &= \frac{\int L(-\frac{1}{\lambda}J(X,U))\nabla_{\theta_j} p(U;\theta)\mathrm{d}U}{\int L(-\frac{1}{\lambda}J(X,U))p(U;\theta)\mathrm{d}U} \\
&= \frac{\int L(-\frac{1}{\lambda}J(X,U))p(U;\theta) \nabla_{\theta_j} \ln p(U;\theta)\mathrm{d}U}{\int L(-\frac{1}{\lambda}J(X,U))p(U;\theta)\mathrm{d}U}
\end{split}
\end{equation}

The term $\nabla_{\theta_j} \ln p(U;\theta)$ can be calculated as:
\begin{equation}
\begin{split}
\nabla_{\theta_j} \ln p(U;\theta) &= \nabla_{\theta_j} \sum_{j=0}^{N-1}\big(\ln h(\mathbf{u}_j) + \theta_j^\mathrm{T}T(\mathbf{u}_j)-\frac{1}{2} \theta_j^\mathrm{T}\theta_j \big) \\
&= T(\mathbf{u}_j) - \theta_j
\end{split}
\end{equation}

Plug it back into the gradient we have:
\begin{equation}
\nabla_{\theta_j} l(\theta) = \frac{\int L(-\frac{1}{\lambda}J(X,U))p(U;\theta) (T(\mathbf{u}_j) - \theta_j) \mathrm{d}U}{\int L(-\frac{1}{\lambda}J(X,U))p(U;\theta)\mathrm{d}U}
\end{equation}

With the gradient the update rule for $\theta$ can be found as:
\begin{equation}
\begin{split}
\theta_j^{k+1} &= \theta_j^k + \\
&\alpha^k\Bigg(\frac{\int L(-\frac{1}{\lambda}J(X,U))p(U;\theta) (T(\mathbf{u}_j) - \theta_j^k) \mathrm{d}U}{\int L(-\frac{1}{\lambda}J(X,U))p(U;\theta)\mathrm{d}U}\Bigg) \\
&= \theta_j^k + \alpha^k\Bigg(\frac{\mathbb{E}\Big[L(-\frac{1}{\lambda}J(X,U)) (T(\mathbf{u}_j) - \theta_j^k) \Big]}{\mathbb{E} \Big[L(-\frac{1}{\lambda}J(X,U))\Big]}\Bigg)
\end{split}
\end{equation}

Then $\bm{\mu}_j$ can be substituted in for $\theta_j$:
\begin{equation}
\begin{split}
&(\bm{\Sigma}_D+I_j\bm{\Sigma}_J)^{-\frac{1}{2}}\bm{\mu}_j^{k+1} = (\bm{\Sigma}_D+I_j\bm{\Sigma}_J)^{-\frac{1}{2}}\bm{\mu}_j^k \\
&+ \alpha^k\Bigg(\frac{\mathbb{E}\Big[L(-\frac{1}{\lambda}J(X,U)) (\bm{\Sigma}_D+I_j\bm{\Sigma}_J)^{-\frac{1}{2}} (\mathbf{u}_j - \bm{\mu}_j^k) \Big]}{\mathbb{E} \Big[L(-\frac{1}{\lambda}J(X,U))\Big]}\Bigg) \\
\end{split}
\end{equation}

The final update law can be obtained as:
\begin{equation}
\begin{split}
\bm{\mu}_j^{k+1} &= \bm{\mu}_j^k + \alpha^k\Bigg(\frac{\mathbb{E}\Big[L(-\frac{1}{\lambda}J(X,U)) (\mathbf{u}_j - \bm{\mu}_j^k) \Big]}{\mathbb{E} \Big[L(-\frac{1}{\lambda}J(X,U))\Big]}\Bigg) \\
&= \bm{\mu}_j^k + \alpha^k\Bigg(\frac{\mathbb{E}\Big[\exp(-\frac{1}{\lambda}J(X,U)) (\epsilon_{D,j}^k+I_j\epsilon_{J,j}^k) \Big]}{\mathbb{E} \Big[\exp(-\frac{1}{\lambda}J(X,U))\Big]}\Bigg)
\end{split}
\end{equation}

The update law (49) is very close to the one (36) obtained from the information theoretic approach. In the case of noise entering the system through control channels only, $\mathbf{B}$ and $\mathbf{H}$ matrices are the same as $\mathbf{G}$, and the matrix transform $\mathcal{G}^{-1}\mathcal{B}$, which maps from state space to control space, goes to identity. Taking step size $\alpha=1$, the two control laws differ only in the extra terms resulted from importance sampling.

With this alternative derivation, convergence of the update law can be shown using the ODE method as described in \cite{zhou2014gradient}.

\begin{algorithm} [!hb]
\caption{MPPI Control on Jump Diffusion}
\begin{algorithmic}
\STATE $\textbf{Given:}$
\STATE $\text{M: Number of samples;}$
\STATE $\text{N: Number of timesteps;}$
\STATE $(\mathbf{u}_0,\mathbf{u}_1,\cdots,\mathbf{u}_{N-1}) \text{: Initial control sequence;}$
\STATE $\mathbf{x}_0 \text{: Initial states;}$
\STATE $\Delta t, \mathbf{f}, \mathbf{G}, \mathbf{B}, \mathbf{H} \text{: System/sampling dynamics;}$
\STATE $\phi, q, \lambda, \mathbf{R} \text{: Cost function parameters;}$
\STATE $c,\bm{\Sigma}_D,\bm{\Sigma}_J, \nu\text{: Noise parameters}$
\STATE $\mathbf{u}_{init} \text{: Value for new control initialization;}$
\WHILE {\textit{task not completed}}
\FOR {$m=0$ \text{to} $M-1$}
\STATE $\text{Update } \mathbf{x}_0 \text{;}$
\STATE $\text{Sample } \mathbf{\varepsilon}^m=\big(\epsilon_0^m,\cdots,\epsilon_{N-1}^m\big), \epsilon_i^m \in \mathcal{N}(0,c \bm{\Sigma}_D);$
\FOR {$i=0$ \text{to} $N-1$}
\STATE $p=U(0,1);$
\IF {$p< \nu \Delta t$}
\STATE $\epsilon_j=\mathcal{N}(0,\bm{\Sigma}_J)$
\STATE $\mathbf{x}_{i+1} = \mathbf{x}_i + (\mathbf{f}+\mathbf{G}\mathbf{u}_i)\Delta t + \mathbf{B} \epsilon_i^m\sqrt{\Delta t} + \mathbf{H}\epsilon_j\sqrt{\Delta t};$
\STATE $\epsilon_i^m = \epsilon_i^m + \epsilon_j;$
\ELSE
\STATE $\mathbf{x}_{i+1} = \mathbf{x}_i + (\mathbf{f}+\mathbf{G}\mathbf{u}_i)\Delta t + \mathbf{B} \epsilon_i^m\sqrt{\Delta t};$
\ENDIF
\STATE $\tilde{S}(X^m) = \tilde{S}(X^m) + \tilde{q}(\mathbf{x}_i,\mathbf{u}_i,\epsilon_i^m);$
\ENDFOR
\ENDFOR
\FOR {$i=0$ to $N-1$}
\STATE $\mathbf{u}_i = \mathbf{u}_i + \frac{\sum_{m=1}^M\exp{\big(-\frac{1}{\lambda}\tilde{S}(X^m)\big)}\frac{\epsilon_i^m}{\sqrt{\Delta t}}}{\sum_{m=1}^M\exp{\big(-\frac{1}{\lambda}\tilde{S}(X^m)\big)}};$
\ENDFOR
\STATE $\text{Execute control policy }\mathbf{u}_0;$
\FOR {$i=0$ to $N-2$}
\STATE $\mathbf{u}_i = \mathbf{u}_{i+1};$
\ENDFOR
\STATE $\mathbf{u}_{N-1} = \mathbf{u}_{init};$
\ENDWHILE
\end{algorithmic}
\end{algorithm}

\section{Model Predictive Control Algorithm}
From both approaches, we get an iterative update law for the optimal control policy at each timestep. This allows for the algorithm to be implemented in a MPC fashion. In the MPC setting, after the optimal control sequence is obtained, only the first control action is executed and re-optimization occurs from the new initial states. Since an entire optimal control sequence is given at every timestep, we can keep the un-executed control sequence to warm start optimization for the next iteration. This is very important for increasing the performance of the algorithm as we are reusing information from previous optimization iterations.

Another key aspect of the algorithm is that its computationally involved parts, namely trajectory propagation and cost computation for each sampled trajectory, can be done in parallel on a GPU. Parallel computation allows us to sample thousands of trajectories at the same time rather than in sequence, with the computation time for each iteration of less than 20 miliseconds in our simulation, which satisfies the real time requirement for a 50 Hz controller. The description of the new Model Predictive Path Integral (MPPI) algorithm is given in Algorithm 1.

The algorithm is based on the assumption that both Gaussian and jump noise affect the states through the control channels. Jump noise is simulated using the zero-one jump law, which states that if $\nu\Delta t \ll 1$, the probability of more than one jump occuring at each timestep can be neglected. A jump timer $p$ is sampled from a uniform distribution to check whether jump occurs at each timestep. When a jump occurs, a zero mean Gaussian vector determines the magnitude of jump noise in each control channel.

\section{Simulation Results}
We compare the MPPI algorithm for jump diffusion processes against the old MPPI algorithm in \cite{williams2016aggressive} that doesn't account for jump noise on a cart pole and quadrotor in simulation with artificial Gaussian and jump noise. To avoid confusion, we refer to the algorithm presented in this paper as the \textcolor{green}{new MPPI} algorithm and the algorithm without jump noise in sampling as the \textcolor{blue}{old MPPI} algorithm. In the trajectory plots the mean trajectory and 95\% confidence interval are plotted. Note that the confidence intervals are not labeled explicitly but shaded with the same color as mean trajectories. The red line indicates the target state.

\subsection{Cart Pole}
We applied the new and old MPPI algorithm on a standard cart pole system in simulation. The task is to swing up and stabilize the cart pole. We used 1000 trajectories during sampling and ran each algorithm for 100 trials. We tested the robustness of both algorithms by varying the jump amplitude and rate while keeping Gaussian noise the same. In Table. \ref{table:1}, we demonstrate the simulation results. The new MPPI algorithm has a higher success rate in stabilizing the cart pole. Specifically, with only small jump noise, both algorithms managed to balance the cart pole. As the jump amplitude increased, both algorithms started to fail, but the new algorithm has a higher success rate of stabilizing the pole than the old algorithm. For a fixed jump amplitude, increasing the jump rate results in lower success rates in both algorithms and vice versa.

\begin{figure}[t]
  \centering
  \includegraphics[width=8.75 cm]{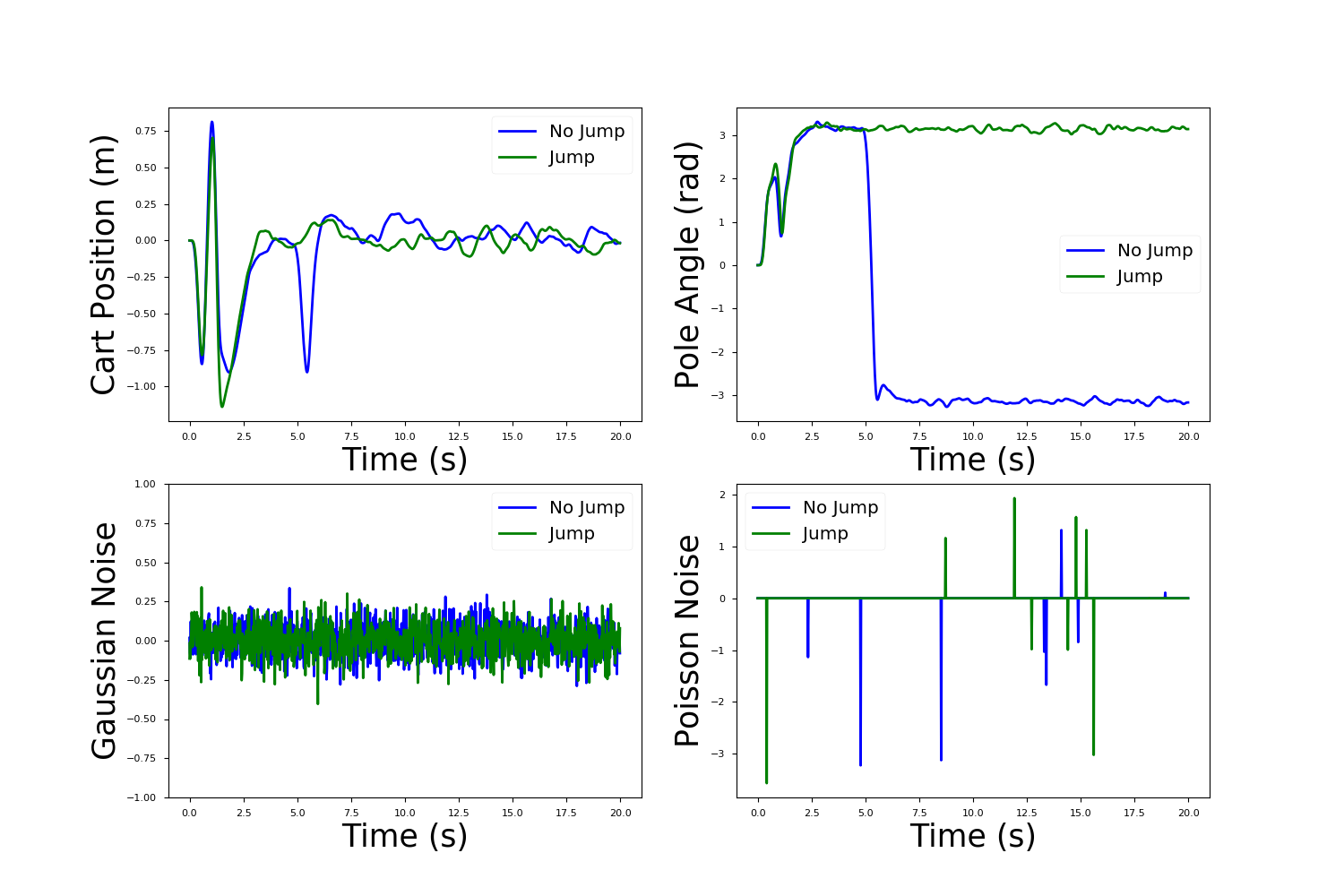}
  \caption{Comparison of \textcolor{blue}{old} and \textcolor{green}{new} algorithm on a cartpole. \textit{Top left: cart position; Top right: pole angle; Botton left: Gaussian noise; Botton right: Poisson noise}.}
  \label{figure:1}
  \vspace{-0.5cm}
\end{figure}

Fig. \ref{figure:1} demonstrates the responses of both algorithms in a trial when the old MPPI algorithm failed. The pole angle plot and the Poisson noise plot show that the old MPPI algorithm failed after a noise spike and had to restabilize the pole. On the other hand, the new algorithm experienced noise spikes of similar magnitude and maintained balance. The cart position plot shows that the new algorithm managed to maintain balance efficiently around the origin.

\begin{table} [!ht]
\hspace{-0.5cm}
\caption{Success rates of \textcolor{blue}{old} and \textcolor{green}{new} algorithm on a cart pole}
\label{table:1}
\begin{center}
\begin{tabular}{l*{2}{c}}
Jump noise & \color{green}New MPPI & \color{blue}Old MPPI \\
\hline
$\nu=0.25,\bm{\Sigma}_J=1$      & \color{green}100\% & \color{blue}100\% \\
$\nu=0.25,\bm{\Sigma}_J=1.5$    & \color{green}96\% & \color{blue}91\% \\
$\nu=0.25,\bm{\Sigma}_J=2$      & \color{green}96\% & \color{blue}81\% \\
$\nu=0.25,\bm{\Sigma}_J=3$     & \color{green}88\%  & \color{blue}61\% \\
$\nu=0.1, \bm{\Sigma}_J=2$    & \color{green}97\% & \color{blue}92\% \\
$\nu=0.5, \bm{\Sigma}_J=2$    & \color{green}91\% & \color{blue}73\%
\end{tabular}
\end{center}
\end{table}

\subsection{Quadrotor}

We also applied both algorithms on a quadrotor system in simulation. The task is to fly from an initial position to a target position. Since it is a more complex system we increased the number of sampling trajectories to 3000 and ran each algorithm for 100 trials. Again we varied the jump amplitude and rate while keeping Gaussian noise the same. Table. \ref{table:2} lists the simulation results. Similar to the cartpole simulation, we found that the new MPPI algorithm has a higher success rate in completing the task. Specifically, with only small jump noise, both algorithms could carry out the task perfectly. As we increased the jump noise amplitude, the failure rate of the old MPPI algorithm increased while the new algorithm maintained perfect task completion rate. Additionally, for a jump amplitude large enough that the old MPPI algorithm has a non zero failure rate, increasing the jump rate further increases the failure rate of the old MPPI algorithm and vice versa.

\begin{table} [t]
\hspace{-0.5cm}
\caption{Success rates of \textcolor{blue}{old} and \textcolor{green}{new} algorithm on a quadrotor}
\label{table:2}
\begin{center}
\begin{tabular}{l*{2}{c}}
Jump noise & \color{green}New MPPI & \color{blue}Old MPPI \\
\hline
$\nu=0.2,\bm{\Sigma}_J=5$      & \color{green}100\% & \color{blue}100\% \\
$\nu=0.2,\bm{\Sigma}_J=10$    & \color{green}100\% & \color{blue}98\% \\
$\nu=0.2,\bm{\Sigma}_J=20$      & \color{green}100\% & \color{blue}97\% \\
$\nu=0.2,\bm{\Sigma}_J=30$     & \color{green}100\%  & \color{blue}87\% \\
$\nu=0.1, \bm{\Sigma}_J=20$    & \color{green}100\% & \color{blue}98\% \\
$\nu=0.5, \bm{\Sigma}_J=20$    & \color{green}100\% & \color{blue}91\%
\end{tabular}
\end{center}
\end{table}

\begin{figure}[t]
  \centering
  \includegraphics[width=8.75 cm]{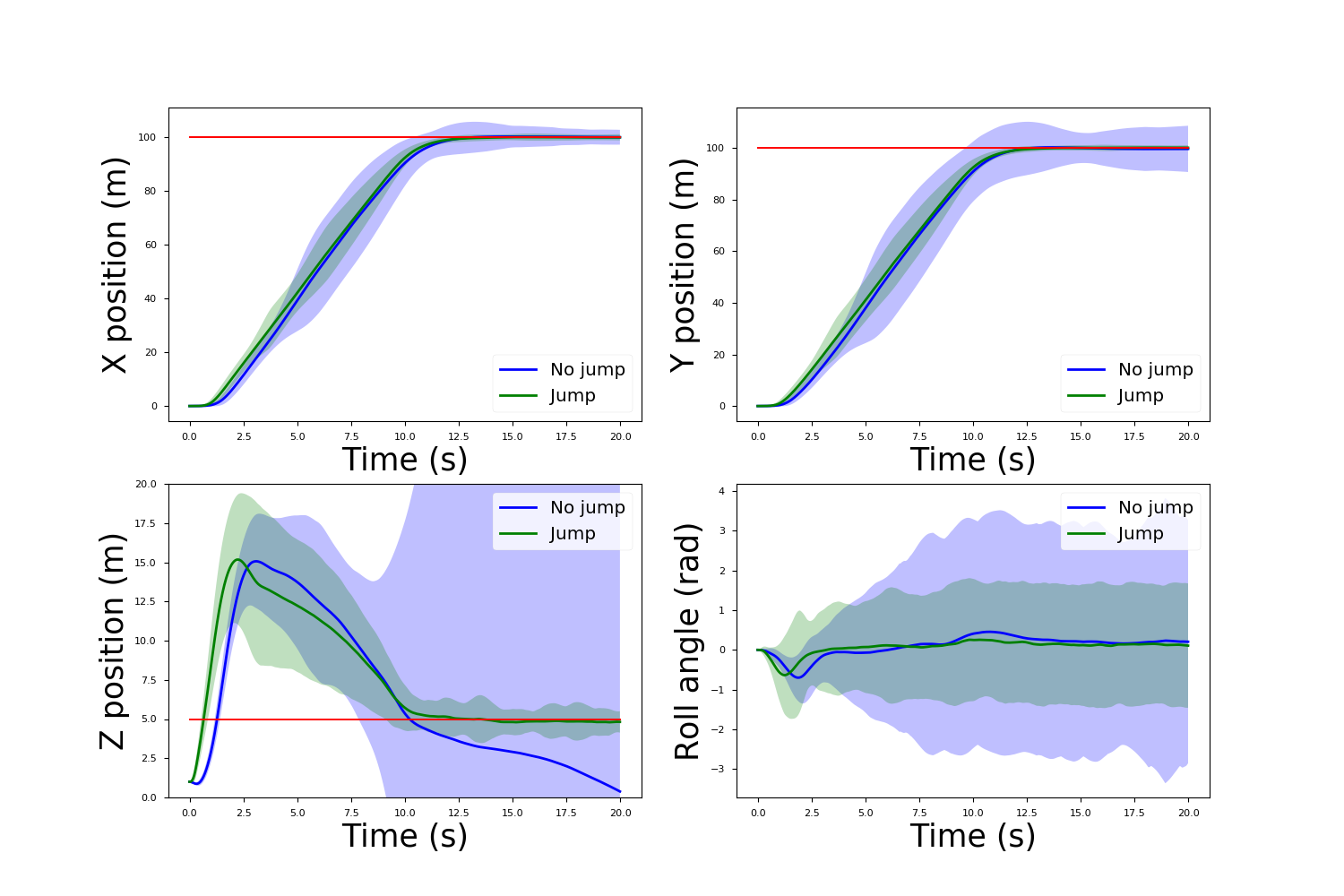}
  \caption{Comparison of \textcolor{blue}{old} and \textcolor{green}{new} algorithm on a quadrotor with 3000 trajectories in sampling. \textit{Top left: x position; Top right: y position; Botton left: z position; Botton right: roll angle}.}
  \label{figure:2}
  \vspace{-0.5cm}
\end{figure}

\begin{figure}[t]
  \centering
  \includegraphics[width=8.75 cm]{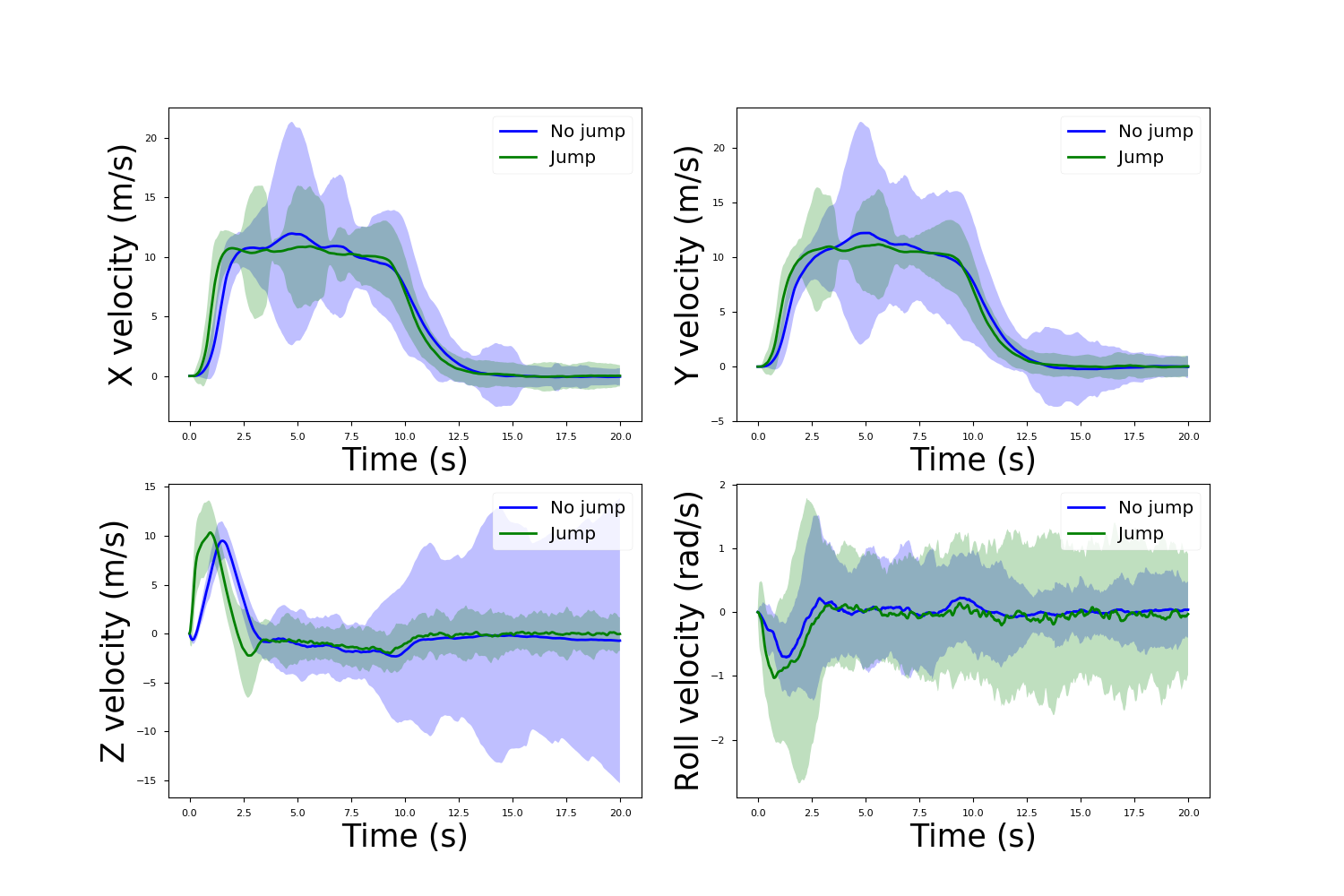}
  \caption{Comparison of \textcolor{blue}{old} and \textcolor{green}{new} algorithm on a quadrotor with 3000 trajectories in sampling. \textit{Top left: x velocity; Top right: y velocity; Botton left: z velocity; Botton right: roll rate}.}
  \label{figure:3}
\end{figure}

In Fig. \ref{figure:2} and \ref{figure:3}, we compare the response of the two algorithms for one test case ($\nu=0.2,\bm{\Sigma}_J=20$). The plots include the mean and 95\% confidence region of the responses. The x and y position plots show that the mean of trajectories resulted from both algorithms follow a similar path to the target, but the variance of trajectories resulted from the old MPPI algorithm is much larger. From the z position plot, we observe that the variance of trajectories resulted from the old MPPI algorithm is significantly larger than the new MPPI algorithm since there are three crash runs. From the x and y velocity plots we find distinct areas where the variance of both algorithms increase. These areas correspond to the high variance regions of the position plots.

\begin{figure}[t]
  \centering
  \includegraphics[width=8.75 cm]{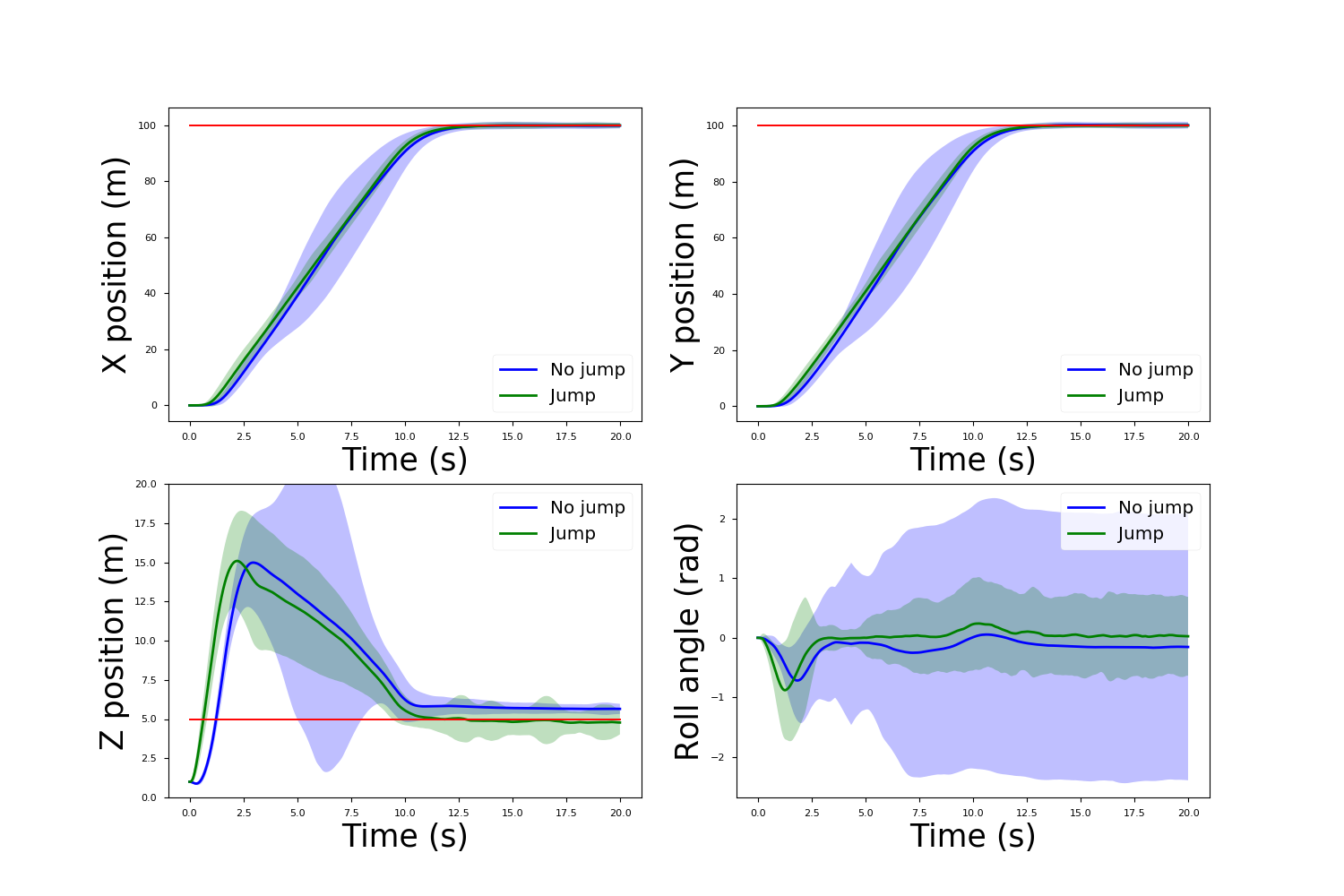}
  \caption{Comparison of \textcolor{blue}{old} and \textcolor{green}{new} algorithm on a quadrotor with 6000 trajectories in sampling. \textit{Top left: x position; Top right: y position; Botton left: z position; Botton right: roll angle}.}
  \label{figure:4}
\end{figure}

\begin{figure}[t]
  \centering
  \includegraphics[width=8.75 cm]{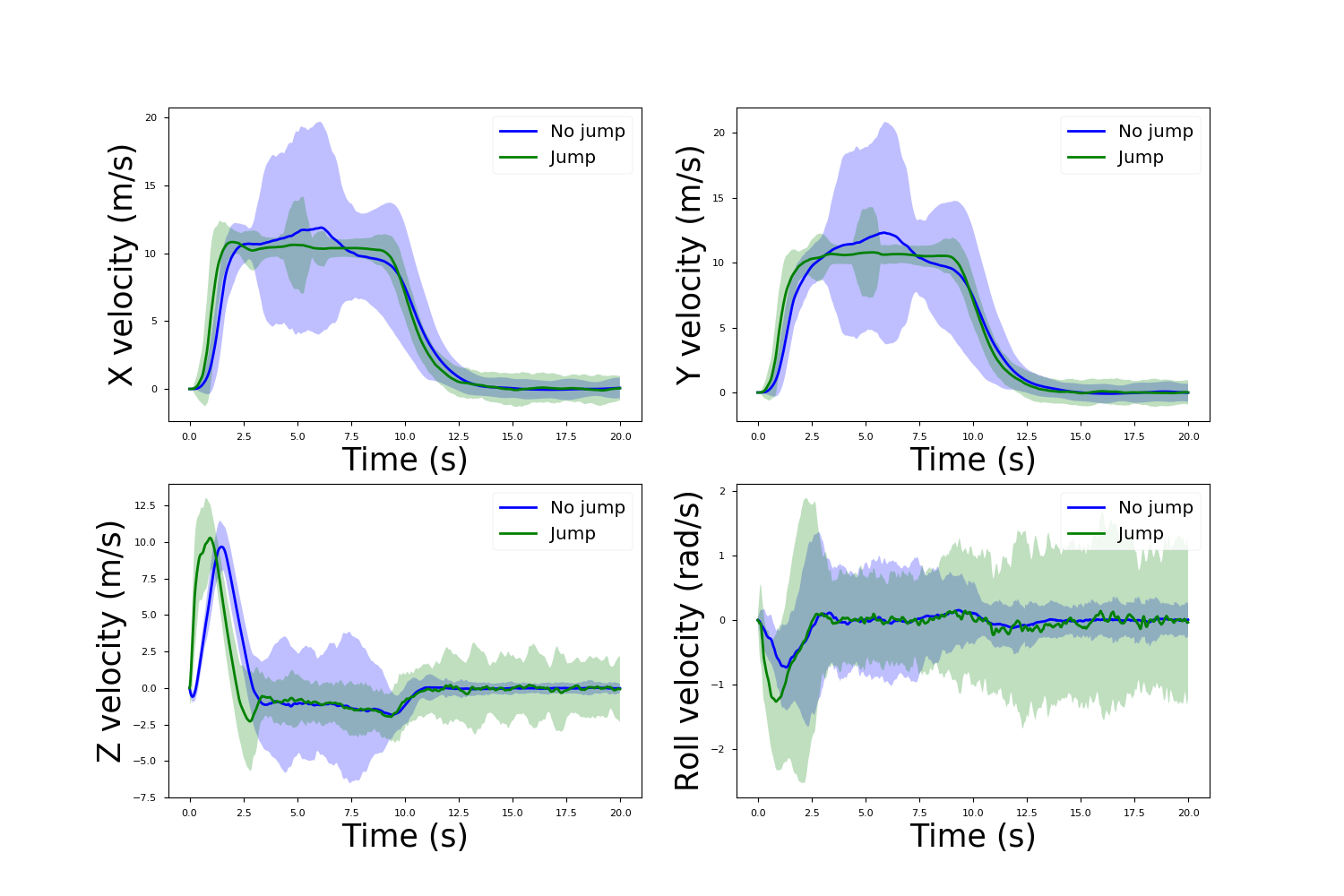}
  \caption{Comparison of \textcolor{blue}{old} and \textcolor{green}{new} algorithm on a quadrotor with 3000 trajectories in sampling. \textit{Top left: x velocity; Top right: y velocity; Botton left: z velocity; Botton right: roll rate}.}
  \label{figure:5}
\end{figure}

We also took two test cases ($\nu=0.2,\bm{\Sigma}_J=5$ and $\nu=0.2,\bm{\Sigma}_J=20$) and ran both algorithms with 6000 sampling trajectories. Fig. \ref{figure:4} and \ref{figure:5} show the response of both algorithms with high jump noise amplitude. Doubling the sampling trajectories resulted in one less crash run for the old MPPI algorithm, while the new MPPI algorithm maintained perfect success rate. From the x and y position plots, we observe that the variance for both algorithms are smaller than the case with fewer sampling trajectories. There is one region in the z position plot where the variance for the old MPPI algorithm increases significantly due to the crash runs. The results suggest that increasing the number of sampling trajectories correspond to a decrease in variance in generated trajectories. The decrease is resulted from better approximation of the expectation with more samples.

\begin{figure}[t]
  \centering
  \includegraphics[width=8.75 cm]{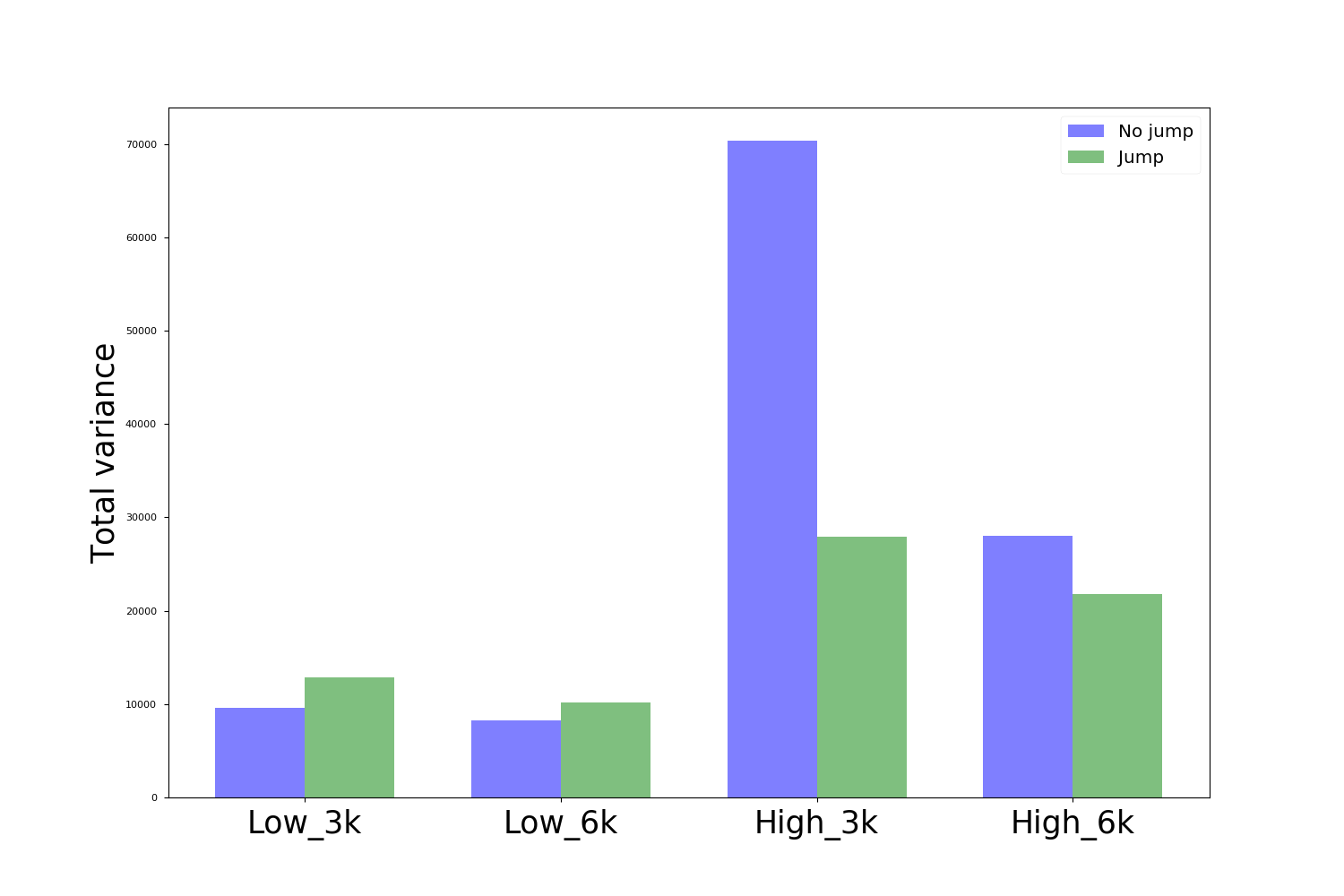}
  \caption{Comparison of total variance of trajectories resulted from \textcolor{blue}{old} and \textcolor{green}{new} algorithm under high ($\bm{\Sigma}_J$=20) and low ($\bm{\Sigma}_J$=5) levels of jump noise amplitude with 3000 and 6000 trajectories in sampling ($\nu$=0.2).}
  \label{figure:6}
\end{figure}

In Fig. \ref{figure:6} we compare the total variance (sum of variance in all states over the entire time horizon for all trajectories) of both algorithms with two jump noise levels using 3000 and 6000 sampling trajectories. We find that with low jump noise amplitude, the old MPPI algorithm results in slightly lower variance than the new MPPI algorithm. The new MPPI algorithm tends to generate trajectories that oscillate around the target location since the dynamics is perturbed more during sampling. For the case of high jump noise amplitude, the difference in variance between the two algorithms is significantly reduced with increased sampling trajectories. This is due to the benefit of better exploration by including jump noise is reduced with increased sampling trajectories.

\section{CONCLUSIONS}

We presented an information theoretic model predictive control algorithm that obtains the optimal control through sampling with Gaussian and compound Poisson noise. We provided an alternative stochastic optimization derivation, from which convergence of the control update law can be proven. We applied the algorithm on cart pole and quadrotor systems with artificially introduced compound Poisson noise and compared its performance to the previously developed algorithm that doesn't include compound Poisson noise in sampling. We demonstrated superior performance of our new algorithm than the old algorithm under large Poisson noise level and comparable performance under low Poisson noise level. Our results suggest that it is important to consider the statistical characteristics of stochastic disturbances in the computation of the optimal control policies.




\addtolength{\textheight}{-12cm}   





\bibliographystyle{unsrt}
\bibliography{references}

\end{document}